\documentclass[11pt,reqno]{amsart}
\usepackage[letterpaper,hmargin=1.05in,vmargin=1in]{geometry}
\usepackage{amssymb,mathtools}
\usepackage[hidelinks]{hyperref}
\hypersetup{pdftitle={Inner derivations and automorphism groups of vertex operator algebras},
 pdfsubject={Derivations, the identity component and finite automorphism groups},
 pdfkeywords={vertex operator algebra, derivation, automorphism group, Zhu algebra}}
\usepackage{microtype}
\newtheorem{theorem}{Theorem}[section]
\newtheorem{proposition}[theorem]{Proposition}
\newtheorem{lemma}[theorem]{Lemma}

\theoremstyle{definition}

\newtheorem{example}[theorem]{Example}
\numberwithin{equation}{section}
\DeclareMathOperator{\Aut}{Aut}
\DeclareMathOperator{\Der}{Der}
\DeclareMathOperator{\Inn}{Inn}
\DeclareMathOperator{\Lie}{Lie}
\DeclareMathOperator{\Tr}{Tr}
\DeclareMathOperator{\End}{End}
\DeclareMathOperator{\Span}{span}
\DeclareMathOperator{\Res}{Res}
\DeclareMathOperator{\wt}{wt}
\newcommand{\C}{\mathbb C}
\newcommand{\Z}{\mathbb Z}
\newcommand{\one}{\mathbf1}
\newcommand{\Lz}{L(0)}
\newcommand{\ct}{\circ_t}
\newcommand{\dtp}{\diamond_t}
\newcommand{\R}{\mathcal R}
\allowdisplaybreaks[2]
\title[Inner derivations and automorphism groups]{Inner derivations and automorphism groups\\ of vertex operator algebras}
\author{Jianzhi Han}
\newcommand{\AuthorAddress}{School of Mathematical Sciences,
  Key Laboratory of Intelligent Computing and Applications
  (Ministry of Education), Tongji University, Shanghai 200092, China}
\newcommand{\AuthorEmail}{jzhan@tongji.edu.cn}
\address{\AuthorAddress}
\email{\AuthorEmail}
\thanks{This work was supported by the National Natural Science Foundation
of China (No.~12271406).}
\date{}
\begin{document}
\begin{abstract}
Let $V$ be a simple, rational, $C_2$-cofinite vertex operator algebra
of CFT type. We prove that the zero mode map $V_1\to\Der(V)$,
$h\mapsto h_0$, is an isomorphism of Lie algebras. We prove surjectivity
using traces on finite sums of weight spaces and Li's module construction,
and injectivity using the vacuum character of Zhu's algebra.
It follows that $\Inn(V)$ is the identity component of $\Aut(V)$ and
has finite index in $\Aut(V)$. Moreover, $\dim\Aut(V)=\dim V_1$, so
$\Aut(V)$ is finite if and only if $V_1=0$.
\end{abstract}
\maketitle

\section{Introduction}

Let $(V,Y,\one,\omega)$ be a vertex operator algebra over $\C$, with vacuum vector $\one$ and conformal vector $\omega$. Write $Y(a,x)=\sum_{n\in\Z}a_nx^{-n-1}$.
A \emph{derivation} of $V$ is a linear map $D:V\to V$ such that $D\one=0$, $D\omega=0$, and
\[
 D(a_n b)=(Da)_n b+a_n(Db)
 \qquad(a,b\in V,\ n\in\Z).
\]
The displayed identity is equivalent to $[D,Y(a,x)]=Y(Da,x)$ for every $a\in V$.  We write $\Der(V)$ for the Lie algebra of derivations, with bracket $[D,E]=DE-ED$.
Taking $a=\omega$ gives $[D,L(n)]=0$ for every $n\in\Z$. In particular, $D$ preserves each weight space $V_n$.

For homogeneous $h\in V_j$, set $o(h)=h_{j-1}$ and extend $o$ linearly to $V$. A derivation $D$ of $V$ is called \emph{inner} if $D=o(h)$ for some $h\in V$. If $V$ is of CFT type, the inner derivations are precisely the operators $u_0$ with $u\in V_1$ \cite[Section~3]{DG}. Each such $u_0$ defines an automorphism $\exp(u_0)$ of $V$. We denote the group of inner automorphisms by
\[
  \Inn(V)=\langle\exp(u_0)\mid u\in V_1\rangle.
\]
Dong and Griess proved that the automorphism group of a finitely generated vertex operator algebra is a complex algebraic group. Let $\Aut(V)_e$ denote the identity component of $\Aut(V)$. In \cite{DG}, it was observed that for a finitely generated rational vertex operator algebra $V$ of CFT type, the equality $\Inn(V)=\Aut(V)_e$ holds if and only if every derivation of $V$ is inner. The authors conjectured that this equality holds for every rational vertex operator algebra of CFT type. It would follow that $\Inn(V)$ has finite index in $\Aut(V)$, although $\Inn(V)$ itself can be infinite, as illustrated by the rank one lattice example at the end of the paper.

\begin{theorem}\label{main}
Let $V$ be a simple, rational, $C_2$-cofinite vertex operator algebra of CFT type. Then
\begin{enumerate}
\item $V_1\longrightarrow\Der(V)$, $h\longmapsto h_0$, is an isomorphism of Lie algebras.
\item $\Aut(V)_e=\Inn(V)$ and $\Aut(V)/\Inn(V)$ is finite.
\item $\dim\Aut(V)=\dim V_1$. In particular, $\Aut(V)$ is finite if and only if $V_1=0$.
\end{enumerate}
\end{theorem}

Put $T=L(-1)$. With Lie bracket $[h,k]=h_0k$ on $V_1$, the map $h\mapsto h_0$ is a Lie algebra homomorphism into $\Der(V)$ by the mode commutator formula and \cite[Section~3]{DG}.

Section~\ref{identity-sec} deduces part (2) of Theorem~\ref{main} from the surjectivity in part (1). Sections 3 and 4 prove surjectivity under the hypotheses of Theorem 1.1.  The first step is to construct a functional on circle products. There are relations among these products, and the main point is to check that the prescribed values respect all of them. We do this in Section~\ref{functional-sec} by taking traces on finite sums of weight spaces. In Section~\ref{realize-sec}, Li's module construction gives a weight one vector that realizes the functional. Section~\ref{group-sec} proves injectivity and gives the criterion for finiteness of the full group.

\section{The identity component}\label{identity-sec}

\begin{proof}[Proof of part (2) of Theorem~\ref{main}]
Rationality implies that $V$ is finitely generated \cite{DZ}. Thus $G=\Aut(V)$ is a complex algebraic group \cite[Theorem~2.4]{DG}. Its Lie algebra acts on $V$ by derivations, and each $D\in\Der(V)$ is the tangent at the identity to $t\mapsto\exp(tD)$, so $\Lie(G)=\Der(V)$ \cite[Section~3]{DG}. By the equivalence in \cite[Introduction]{DG}, $G_e=\Inn(V)$ if every derivation has the form $h_0$ for some $h\in V_1$. Section~\ref{realize-sec} proves precisely this surjectivity. The quotient $G/\Inn(V)=G/G_e$ is finite because a complex algebraic group has finitely many connected components.
\end{proof}

\section{A functional on circle products}\label{functional-sec}

Throughout Sections~\ref{functional-sec} and \ref{realize-sec}, $V$ satisfies the hypotheses of Theorem~\ref{main}. Thus $V=\bigoplus_{n\ge0}V_n$, $V_0=\C\one$, every admissible $V$-module is completely reducible, and $V/C_2(V)$ is finite dimensional. Fix $D\in\Der(V)$ and put $H_t=\Lz+tD$ for $t\in\C$. For $v=\sum_{n\ge0}v^{(n)}$ with $v^{(n)}\in V_n$, define $\chi:V\to\C$ by $v^{(0)}=\chi(v)\one$. Define
\begin{align*}
 a\ct b&=\Res_x x^{-2}Y((1+x)^{H_t}a,x)b,\\
 a\dtp b&=\Res_x (1+x)x^{-2}Y((1+x)^{H_t}a,x)b.
\end{align*}
For $t=0$, we write $\circ$ and $\diamond$ for $\circ_0$ and $\diamond_0$, respectively. In particular, $O(V)=\Span\{a\circ b\}$ \cite{Zhu}. For $a\in V_m$, $(1+x)^{H_t}a=(1+x)^m\exp(t\log(1+x)D)a$. The coefficient of $t^j$ is divisible by $x^j$. If $b\in V_n$, truncation therefore makes each residue a polynomial in $t$ of degree at most $m+n+1$. Define
\[
 C_D(a,b)=\left.\frac{d}{dt}(a\ct b)\right|_{t=0},
 \qquad \beta_D(a,b)=\chi(C_D(a,b)).
\]
By \cite[Theorem~2.1.2]{Zhu}, $A(V)$ acts on $V_0=\C\one$ by zero modes. The creation formula gives $o(a)\one=\chi(a)\one$, so $\chi$ is the corresponding character and kills $O(V)$. For $a\in V_m$ and $b\in V$, Zhu's product is
\[
 a*b=\Res_x x^{-1}(1+x)^mY(a,x)b.
\]
Thus $a\diamond b=a\circ b+a*b$, and $\chi(a\diamond b)=\chi(a)\chi(b)$. Also $\chi H_t=0$ since $H_t$ preserves weights and kills $\one$.

\begin{proposition}\label{functional}
There is a linear functional $\lambda:V\to\C$ such that
\begin{equation}\label{lambda}
  \lambda(a\circ b)=-\beta_D(a,b)\qquad(a,b\in V).
\end{equation}
\end{proposition}

To construct $\lambda$, we use Zhu's notation $Y[a,z]=Y(e^{z\Lz}a,e^z-1)=\sum_j a[j]z^{-j-1}$. Let $E_{2j}$ denote the normalized Eisenstein series with constant term one. For $j\ge2$, $E_{2j}\in\C[E_4,E_6]$, and $E_2,E_4,E_6$ are algebraically independent \cite{KZ,Zagier}. In particular, if $A,B\in\C(X,Y)$ and $A(E_4,E_6)=E_2B(E_4,E_6)$, then $A=B=0$.

Let $q$ be the formal variable in the Fourier expansions of the
Eisenstein series. Write $w=e^z$ and set
\begin{align*}
 k(z)&=\frac12+\frac{w}{1-w}
  +\sum_{d\ge1}\left(\frac{q^dw}{1-q^dw}-\frac{q^d/w}{1-q^d/w}\right),\\
 P(z)&=k'(z)+E_2/12,\qquad K(z)=k(z)+E_2z/12 .
\end{align*}
The prime denotes $d/dz$. Expanding the geometric series and using the Bernoulli numbers $B_{2j}$ gives
\[
 k(z)=-z^{-1}-\sum_{j\ge1}\frac{B_{2j}}{(2j)!}E_{2j}z^{2j-1}.
\]
Thus the leading terms of $P,K$ are $z^{-2},-z^{-1}$, respectively. Their other coefficients are polynomials in $E_4,E_6$. Expanding each coefficient of $q^d$ in $k$ at $w=0$ and at $w=\infty$ gives the following series $k^{\rm in}$ and $k^{\rm out}$, respectively.
\[
 k^{\rm in}(w)=\frac12+\sum_{j\ge1}\frac{w^j-q^jw^{-j}}{1-q^j},
 \qquad
 k^{\rm out}(w)=-\frac12+\sum_{j\ge1}\frac{q^jw^j-w^{-j}}{1-q^j}.
\]
In particular $k^{\rm in}(qw)=k^{\rm out}(w)+1$ and $P^{\rm in}(qw)=P^{\rm out}(w)$.

For $f=1$ or $P$ and $a,b\in V$, write $\alpha=(f,a,b)$ and put
\[
 r_\alpha=\Res_z f(z)Y[a,z]b,\qquad
 s_{p,\alpha}=\Res_z f(z)p(z)Y[Da,z]b\quad(p=k,K,z).
\]
All these residues lie in $U[[q]]$ for some finite dimensional $U\subset V$ depending on $a,b$. Indeed, $f(z)$ and $f(z)p(z)$ have poles of order at most three
at $z=0$. Thus only terms $a[j]b$ and $(Da)[j]b$ with $j\ge -3$ can contribute to the
residues, and both vanish for all
sufficiently large $j$.  For $r_\alpha,s_{K,\alpha},s_{z,\alpha}$ we can also regard the residues as elements of $V\otimes\C[X,Y]$ by writing $X,Y$ for $E_4,E_6$.

We study the relations among the residues $r_\alpha$ by taking traces on irreducible ordinary $V$-modules. Fix such a module $W$ and write $W(n)$ for its $n$th level, with $W(0)$ denoting the lowest weight space.
Set $E=W\oplus W$ and define
\[
 Y_E(a,x)=
 \begin{pmatrix}
  Y_W(a,x) & Y_W(Da,x)\\
  0        & Y_W(a,x)
 \end{pmatrix}
 \qquad (a\in V).
\]
Since $D\one=0$, the vacuum acts as the identity on $E$.
The diagonal entries satisfy the Jacobi identity on $W$. For the upper
right entry, add the Jacobi identities on $W$ for the pairs $(Da,b)$
and $(a,Db)$ and use $Y(Da,x)b+Y(a,x)Db=D(Y(a,x)b)$. Since $D$ preserves weights,
$E(n)=W(n)\oplus W(n)$ gives an admissible grading. Thus we have
a short exact sequence of $V$-modules
\[
 0\longrightarrow W\xrightarrow{\iota}E\xrightarrow{\pi}W
 \longrightarrow 0,
 \qquad \iota(v)=(v,0),\quad \pi(u,v)=v.
\]
By rationality, $\pi$ has a $V$-module section $s:W\to E$.
Since $\pi s=\mathrm{id}_W$, we can write $s(v)=(\mathcal Bv,v)$ for a
linear operator $\mathcal B$ on $W$.
The first components of $s(Y_W(a,x)v)=Y_E(a,x)s(v)$ give
\[
 \mathcal B Y_W(a,x)v=Y_W(a,x)\mathcal Bv+Y_W(Da,x)v.
\]
Thus
\begin{equation}\label{B}
 [\mathcal B,Y_W(a,x)]=Y_W(Da,x).
\end{equation}
Taking $a=\omega$ and using $D\omega=0$, we obtain $[\mathcal B,L(0)]=0$.
Hence $\mathcal B$ preserves each $W(n)$.

For a homogeneous $a$ put $A_j(a)=a^W_{\wt a-1+j}$ and extend linearly. Thus $A_j(a)$ has degree $-j$, and $o(a)=A_0(a)$. For a series $v(q)=\sum_{d\ge0}v_dq^d$ define
\[
 \mathcal T_{W,N}(v)=\sum_{n=0}^N\Tr_{W(n)}o(v_{N-n}).
\]
\begin{lemma}\label{trace-lemma}
Let $W$ be an irreducible ordinary $V$-module and let $N\ge0$ be an integer.
For every $\alpha=(f,a,b)$ with $f\in\{1,P\}$ and $a,b\in V$, we have
\begin{equation}\label{trace1}
 \mathcal T_{W,N}(r_\alpha)=0.
\end{equation}
Moreover, for any finite family of such triples and coefficients
$c_\alpha(q)\in\C[[q]]$, we have
\begin{equation}\label{trace2}
 \sum_\alpha c_\alpha(q)r_\alpha=0
 \quad\Longrightarrow\quad
 \mathcal T_{W,N}\!\left(\sum_\alpha c_\alpha(q)s_{k,\alpha}\right)=0.
\end{equation}
\end{lemma}

\begin{proof}
Equation~\eqref{trace1} is also a consequence of
\cite[(4.3.14) and~(4.3.16)]{Zhu}. We derive both identities
from the same graded trace calculation.

For $h\in\C[w,w^{-1},(w-1)^{-1}]$, let $h_j^{\rm out}$ and
$h_j^{\rm in}$ be the coefficients of $w^j$ in the Laurent
expansions of $h$ at $w=\infty$ and $w=0$, respectively. Then
\begin{equation}\label{jac}
 o\!\left(\Res_z h(e^z)Y[a,z]b\right)
 =\sum_{j\in\Z}\left(
     h_j^{\rm out}A_j(a)A_{-j}(b)
     -h_j^{\rm in}A_{-j}(b)A_j(a)
   \right).
\end{equation}
Indeed, for $a\in V_p$, $b\in V_r$, and $h(w)=w^\ell(w-1)^s$,
this is \cite[(1.2.9)]{Zhu} with
$(m,l,n)=(p+\ell-1,s,r-\ell-s-1)$. Since these functions span
$\C[w,w^{-1},(w-1)^{-1}]$, the identity follows by linearity.
On each $W(n)$, only finitely many terms contribute.

For an operator $X$ preserving each $W(n)$, put
\[
 \tau(X)=\sum_{n\ge0}q^n\Tr_{W(n)}X,
\]
and extend $\tau$ to series of such operators by $\C[[q]]$-linearity,
taking each coefficient separately. In particular,
\[
 [q^N]\tau(o(v(q)))=\mathcal T_{W,N}(v(q))
 \qquad (v(q)\in V[[q]]).
\]
We use the following form of cyclicity. If $U$ has degree $-j$
and $V'$ has degree $j$, then
\[
 \tau(V'U)=q^j\tau(UV').
\]
Indeed, $U$ maps $W(n)$ to $W(n-j)$, so
$\Tr_{W(n)}(V'U)=\Tr_{W(n-j)}(UV')$, with $W(m)=0$ for $m<0$.
Multiplying by $q^n$ and summing over $n$ gives the formula.

Fix $\alpha=(f,a,b)$ and put
\[
 U_j=A_j(a),\qquad V_j=A_{-j}(b),\qquad g=fk.
\]
Use the same notation for the Laurent coefficients of $f$ and $g$,
now taken coefficientwise in $q$. Since
$f^{\rm in}(qw)=f^{\rm out}(w)$ and
$k^{\rm in}(qw)=k^{\rm out}(w)+1$, comparison of the coefficients
of $w^j$ gives
\[
 q^j f_j^{\rm in}=f_j^{\rm out},\qquad
 q^j g_j^{\rm in}=g_j^{\rm out}+f_j^{\rm out}.
\]

Let $C$ be an operator preserving each $W(n)$. Apply \eqref{jac}
with $h=f$ and take the trace after multiplying on the left by $C$.
Since $\tau(CV_jU_j)=q^j\tau(U_jCV_j)$, we get
\begin{align*}
 \tau(C\,o(r_\alpha))
 &=\sum_j\bigl(f_j^{\rm out}\tau(CU_jV_j)
              -q^j f_j^{\rm in}\tau(U_jCV_j)\bigr)\\
 &=\sum_j f_j^{\rm out}\tau([C,U_j]V_j).
\end{align*}
Taking $C=1$ gives \eqref{trace1}. For $C=\mathcal B$,
comparing modes in \eqref{B} gives $[\mathcal B,U_j]=A_j(Da)$.
Apply \eqref{jac} to $g=fk$ and the pair $Da,b$. Cyclicity gives
\[
 \tau(o(s_{k,\alpha}))
 =\sum_j(g_j^{\rm out}-q^j g_j^{\rm in})\tau(A_j(Da)V_j)
 =-\tau(\mathcal B\,o(r_\alpha)).
\]
These sums are finite at each power of $q$. For $j>0$, the outer
coefficients are divisible by $q^j$. For $j<0$, $V_j$ vanishes
on $W(n)$ if $n<-j$. Multiplying the last identity by
$c_\alpha(q)$ and summing over $\alpha$, the right hand side
vanishes when $\sum_\alpha c_\alpha(q)r_\alpha=0$.
Taking the coefficient of $q^N$ gives \eqref{trace2}.
\end{proof}

We now relate the residues $r_\alpha$ to Zhu's circle products by setting $q=0$. Put $A=\C[X,Y]$ and let $\R_A$ be the $A$-span of the $r_\alpha$ in $V\otimes A$. A subscript $*$ will denote evaluation at $(X,Y)=(1,1)$, equivalently at $q=0$. In particular,
\[
 k_*'(z)=\frac{e^z}{(e^z-1)^2},\qquad P_*=k_*'+1/12 .
\]
At $q=0$, the two choices $f=1$ and $f=P$ give
\[
 (r_{(1,a,b)})_* = a[0]b,\qquad
 (r_{(P,a,b)})_* = a\circ b+\frac{1}{12}a[0]b.
\]
Indeed, $P_*=k_*'+1/12$, and the substitution $x=e^z-1$ gives
$\Res_z k_*'(z)Y[a,z]b=a\circ b$. Thus
\[
 \R_*:=\Span\{(r_\alpha)_*\}
      =O(V)+\Span\{a[0]b\mid a,b\in V\}.
\]
By \cite[(2.1.4)]{Zhu}, the image of $a[0]b$ in $A(V)$ is
$\bar a*\bar b-\bar b*\bar a$. Hence
\begin{equation}\label{cusp}
 V/\R_*\simeq A(V)/[A(V),A(V)]
\end{equation}
as vector spaces. Rationality and ordinary Zhu theory give
$A(V)\simeq\bigoplus_{i=1}^r\End W_i(0)$, where $W_1=V,W_2,\ldots,W_r$ are the inequivalent irreducible modules \cite{DLM,Zhu}. Choose $b_i\in V$ so that $o(b_i)$ acts as the identity on $W_i(0)$ and as zero on $W_j(0)$ for $j\ne i$. They form a basis modulo $\R_*$ and
\begin{equation}\label{tops}
 \Tr_{W_i(0)}o(b_j)=\delta_{ij}\dim W_i(0).
\end{equation}
Write $B_0=\Span\{b_1,\ldots,b_r\}$.

Identify $A=\C[X,Y]$ with $\C[E_4,E_6]$ and set $Q_A=(V\otimes A)/\R_A$. For $g=h=1$, the defining relations of $O(g,h)$ in \cite[(5.1)--(5.2)]{DLMT} are precisely the residues with kernels $1$ and $P$. Thus $Q_A$ is finitely generated over $A$ by \cite[Lemma~5.2]{DLMT}. Put $\mathfrak m=(X-1,Y-1)$, $S=A_{\mathfrak m}$ and $Q_S=Q_A\otimes_A S$. Since $Q_S/\mathfrak m Q_S\simeq V/\R_*$, Nakayama's lemma and \eqref{cusp} show that the classes of $b_1,\ldots,b_r$ generate $Q_S$.

For each irreducible module $W_i$, let $\tau_i$ denote the graded
trace on $W_i$. For $v(q)=\sum_{d\ge0}v_dq^d\in V[[q]]$, we have
\[
 \tau_i(o(v(q)))
 =\sum_{n,d\ge0}q^{n+d}\Tr_{W_i(n)}o(v_d).
\]
We claim that if $v(q)\in B_0[[q]]$ and $\tau_i(o(v(q)))=0$ for
every $i$, then $v(q)=0$. Indeed, suppose that $v_N$ is the first
nonzero coefficient of $v(q)$. Write $v_N=\sum_j c_jb_j$ and choose
$i$ with $c_i\ne0$. Since $v_d=0$ for $d<N$, only the trace on
$W_i(0)$ contributes to the coefficient of $q^N$. Thus \eqref{tops}
gives
\[
 [q^N]\tau_i(o(v(q)))
 =\Tr_{W_i(0)}o(v_N)
 =c_i\dim W_i(0)\ne0,
\]
which proves the claim.

Substitution of $X=E_4$ and $Y=E_6$ embeds $S$ in $\C[[q]]$,
since $E_4,E_6$ are algebraically independent and every denominator
in $S$ has nonzero constant term after substitution. We may therefore
regard $B_0\otimes S$ as a subspace of $B_0[[q]]$.
Let $\R_S$ be the $S$-span of the $r_\alpha$. By \eqref{trace1}
and $\C[[q]]$-linearity, $\tau_i(o(r))=0$ for every $r\in\R_S$
and every $i$. The claim now gives $(B_0\otimes S)\cap\R_S=0$.
Since the classes of $b_1,\ldots,b_r$ generate $Q_S$, we obtain
\[
 V\otimes S=(B_0\otimes S)\oplus\R_S.
\]
Let $\pi:V\otimes S\to B_0\otimes S$ be the projection onto the
first summand.

\begin{proof}[Proof of Proposition~\ref{functional}]
Extend $\chi$ $S$-linearly. On the spanning vectors of $V\otimes S$ prescribe
\begin{equation}\label{ell}
 \ell(b_i)=0,\qquad \ell(r_\alpha)=-\chi\pi(s_{z,\alpha}).
\end{equation}
Since $V\otimes S=(B_0\otimes S)\oplus\R_S$, it suffices to check
the relations among the $r_\alpha$. Given a relation
$\sum_\alpha f_\alpha r_\alpha=0$ over $S$, put
$x_p=\sum_\alpha f_\alpha s_{p,\alpha}$ for $p=k,K,z$. Then
$x_k=x_K-E_2x_z/12$, and \eqref{trace2} gives
$\tau_i(o(x_k))=0$ for every $i$. By \eqref{trace1}, replacing
$x_K$ and $x_z$ by their projections does not change their traces.
Hence every $\tau_i$ vanishes on the zero mode of
$\pi(x_K)-E_2\pi(x_z)/12$. The traces separate $B_0[[q]]$, so
\[
 \pi(x_K)=\tfrac{1}{12}E_2\pi(x_z).
\]
Both projections have coordinates in $S$ in the basis $b_1,\ldots,b_r$.
Since $E_2\notin\operatorname{Frac}(S)$, each coordinate of
$\pi(x_z)$ is zero. Thus \eqref{ell} defines an $S$-linear functional.

Set $\lambda(v)=\ell(v)_*$ for $v\in V$. By \eqref{cusp}, $\chi$ kills $\R_*$, hence $\chi\pi_*=\chi$. Specializing the two defining relations for $f=P$ and $f=1$, and subtracting $1/12$ times the second from the first, gives
\[
 \lambda(a\circ b)=-\chi\!\left(\Res_z zk_*'(z)Y[Da,z]b\right).
\]
Since $z=\log(1+x)$, the residue is $C_D(a,b)$. This proves \eqref{lambda}.
\end{proof}

\section{Realizing the functional by a weight one vector}\label{realize-sec}

Put $S_t=T+H_t$, $S_0=T+\Lz$, and define
\[
 Y_t(a,z)=Y[e^{ztD}a,z],\qquad
 \phi(z)=\frac{e^z}{e^z-1},\qquad \gamma(z)=\phi(z)^2-\phi(z)=-\phi'(z).
\]
Regard $t$ as an indeterminate and extend $Y_t$ by $\C[t]$-bilinearity
to $V[t]=V\otimes_{\C}\C[t]$. For $a\in V_p$ and $b\in V_r$,
$z^{p+r}Y[D^ka,z]b\in V[[z]]$ for every $k\ge0$, since $D^ka\in V_p$.
Thus each coefficient of
\[
 Y_t(a,z)b=\sum_{k\ge0}\frac{t^kz^k}{k!}Y[D^ka,z]b
\]
is a polynomial in $t$, and $Y_t(a,z)b\in V[t]((z))$.
The $D$-orbits of $a$ and $b$ span finite dimensional subspaces, so
locality for their bases gives locality for $Y_t$.
Since $D$ is a derivation of the bracket vertex algebra and commutes
with $S_0$, we have $Y_t(a,z)\one=e^{zS_t}a$ and
$[S_t,Y_t(a,z)]=\partial_zY_t(a,z)=Y_t(S_ta,z)$.
Together with $Y_t(\one,z)=1$, these identities give a vertex
algebra on $V[t]$. Changing variables gives
\[
 a\ct b=\Res_z\gamma(z)Y_t(a,z)b,\qquad
 a\dtp b=\Res_z\phi(z)^2Y_t(a,z)b,\qquad S_ta=a\ct\one.
\]

Let $O_t=\Span_{\C[t]}\{a\ct b\mid a,b\in V[t]\}$. By bilinearity,
any polynomial coefficient can be absorbed into the first entry, so
$O_t$ consists of finite sums of circle products in $V[t]$.

\begin{lemma}\label{products}
For every $a\in V[t]$, we have $a\dtp O_t\subseteq O_t$.
Moreover, for all $a,b\in V[t]$,
\begin{align}
 (S_ta)\dtp b&=((S_t+1)a)\ct b,\label{ST}\\
 H_t(a\dtp b)&=((H_t+1)a)\ct b+a\dtp H_tb.\label{HT}
\end{align}
\end{lemma}

\begin{proof}
Put
\[
 a*_t b=\Res_z\phi(z)Y_t(a,z)b.
\]
For $f=\phi$ and $g=\phi'=-\gamma$, the span of the $g$-products
in $V[t]$ is $O_t$. By the associative quotient construction of
\cite[Section~3]{vEH} with $c=1$, $*_t$ induces a product on
$V[t]/O_t$. Thus $a*_t O_t\subseteq O_t$. Since
$a\dtp b=a*_t b+a\ct b$, the same inclusion holds for $\diamond_t$.

Since $\phi^2-\gamma=\phi$ and the residue of a formal derivative
is zero, we have
\[
 (S_ta)\dtp b-(S_ta)\ct b
 =\Res_z\phi(z)\partial_zY_t(a,z)b
 =-\Res_z\phi'(z)Y_t(a,z)b=a\ct b.
\]
This proves \eqref{ST}. For \eqref{HT}, observe that
$[T,H_t]=-T$ gives $T e^{zH_t}=e^{-z}e^{zH_t}T$, and hence
\[
 [T,Y_t(a,z)]=e^{-z}Y_t(Ta,z).
\]
Since $e^{-z}\phi^2=\gamma$, it follows that
$T(a\dtp b)-a\dtp Tb=(Ta)\ct b$. The translation operator $S_t$
satisfies $S_t(a\dtp b)-a\dtp S_tb=(S_ta)\dtp b$. Subtracting
and using \eqref{ST} gives
\[
 H_t(a\dtp b)-a\dtp H_tb
 =(S_ta)\dtp b-(Ta)\ct b=((H_t+1)a)\ct b.
\]
\end{proof}

Let
\[
 N=\bigoplus_{n\ge0}N(n),\qquad N(n)=A^\diamond_{n,1}(V),
\]
be the admissible $V$-module of \cite[Theorem~4.8]{Li}.
By \cite[Definitions~2.10 and~2.12, (2.23)]{Li}, we have $N(0)=V/J$,
where $J=O'_{0,1}(V)$ is spanned by $(S_0+1)a$ and the residues
$\Res_z\phi(z)\gamma(z)Y_0(a,z)b$. Since $(\phi^2)'=-2\phi\gamma$
and a formal derivative has zero residue,
\[
 2\Res_z\phi(z)\gamma(z)Y_0(a,z)b
 =\Res_z\phi(z)^2\partial_zY_0(a,z)b=(S_0a)\diamond b .
\]
Thus
\[
 J=\Span\{(S_0+1)a,\ (S_0a)\diamond b\mid a,b\in V\}.
\]

Write $[a]=a+J\in N(0)$ and let
$w=\one+O'_{1,1}(V)\in N(1)$. For homogeneous $a$, taking $m=1$
and $n=0$ in the proof of \cite[Lemma~4.10]{Li} gives
\[
 a^N_{\wt(a)}\bigl(\one+O'_{1,1}(V)\bigr)=a+O'_{0,1}(V).
\]
On $N(0)$, the left action in \cite[Definition~3.12]{Li} is induced
by the product $\diamond$. Hence
\begin{equation}\label{li-action}
 a_m^N w=[a]\quad(a\in V_m),\qquad
 o_N(b)[a]=[b\diamond a]\quad(a,b\in V).
\end{equation}

\begin{lemma}\label{realize}
Let $\lambda$ be the functional in Proposition~\ref{functional}.
There exists $u\in V_1$ such that
\[
 \lambda(v)=\chi(u[1]v)\qquad (v\in O(V)).
\]
\end{lemma}

\begin{proof}
For $y=y(t)\in V[t]$, put
\[
 \rho(y)=\lambda(y(0))+\chi(y'(0)).
\]
This functional satisfies
\begin{equation}\label{dif}
 \rho(O_t)=0,\qquad \rho(H_t y)=\lambda(\Lz y(0)).
\end{equation}
For the first, differentiate $a(t)\ct b(t)$ at $t=0$. The terms
from $a'(0)$ and $b'(0)$ lie in $O(V)$ and are killed by $\chi$.
The remaining term has vacuum coefficient
$\beta_D(a(0),b(0))=-\lambda(a(0)\circ b(0))$ by \eqref{lambda}.
The second property follows from $\chi H_t=0$.

If $y(0)\in O(V)$, write $y=r+tv$ with $r\in O_t$ and $v\in V[t]$.
Such an $r$ is obtained by writing $y(0)$ as a sum of circle products
and replacing $\circ$ by $\circ_t$. For $a\in V$, Lemma~\ref{products}
gives $a\dtp r\in O_t$. Thus $(a\dtp y)(0)\in O(V)$ and
\[
 \rho(a\dtp y)=\chi(a\diamond v(0))
 =\chi(a)\chi(v(0))=\chi(a)\rho(y).
\]

Set $U=\{[b]\in N(0)\mid \Lz b\in O(V)\}$. To see that this
condition is independent of the representative, lift the generators
of $J$ to $j(t)=(S_t+1)a$ and $j(t)=(S_ta)\dtp b$. In both cases
$H_tj(t)\in O_t$. For the first, $[H_t,S_t]=T$ gives
$H_t(S_t+1)a=S_t(H_t+1)a\in O_t$. For the second, \eqref{HT} gives
$H_tj(t)\equiv(S_ta)\dtp H_tb\equiv0\pmod{O_t}$ by \eqref{ST}.
Evaluating at zero gives $\Lz J\subseteq O(V)$, and applying
\eqref{dif} gives $\lambda(\Lz J)=0$. Hence
$\eta_U([b])=\lambda(\Lz b)\one$ is well defined on $U$.

Now let $[b]\in U$. By \eqref{HT},
$H_t(a\dtp b)\equiv a\dtp H_tb\pmod{O_t}$. Since
$(H_tb)(0)=\Lz b\in O(V)$, the preceding product calculation gives
$\Lz(a\diamond b)\in O(V)$ and
\[
 \lambda(\Lz(a\diamond b))=\rho(a\dtp H_tb)
 =\chi(a)\rho(H_tb)=\chi(a)\lambda(\Lz b).
\]
By \eqref{li-action}, $U$ is an $A(V)$-submodule and $\eta_U$ is
$A(V)$-linear.

Let $\Omega(N)$ be the subspace annihilated by all modes of negative
degree. We have $U\subseteq N(0)\subseteq\Omega(N)$.
Zero modes preserve degree, so semisimplicity of $A(V)$ gives a
graded $A(V)$-module complement of $U$ in $\Omega(N)$.
Extend $\eta_U$ by zero on this complement.
The resulting $A(V)$-module map $\eta:\Omega(N)\to V_0$
satisfies
\begin{equation}\label{projection}
 \eta([b])=\lambda(\Lz b)\one
 \qquad\text{whenever }\Lz b\in O(V).
\end{equation}
Since $V$ is simple, \cite[Proposition~3.4(i)]{DLM} gives
$\Omega(V)=V_0$. By rationality, the equivalence in
\cite[Theorem~4.10(b), $n=0$]{DLM} lifts $\eta$ to a
$V$-module homomorphism $p:N\to V$.

Let $u$ be the weight one component of $p(w)$. For $a\in V_m$,
\eqref{li-action} gives $a_m p(w)=p([a])=\eta([a])\in V_0$.
Since $a_m$ lowers weight by one, taking the weight zero component
gives $a_m u=\eta([a])$.

In any vertex operator algebra of CFT type, if $h\in V_1$ and
$a\in V_m$ with $m>0$, then
\begin{equation}\label{weight-one-pairing}
 \chi(h[1]a)\one
 =\frac{(-1)^{m-1}}{m}\chi(h_m a)\one
 =\frac1m a_m h .
\end{equation}
The first equality follows from $h[1]=\Res_x\log(1+x)Y(h,x)$,
since only $h_m a$ has weight zero. The second follows from
skew symmetry and $\chi T=0$.

Now let $v\in O(V)$. Since $\chi(v)=0$, we can write
$b=\Lz^{-1}v=\sum_{m>0}b^{(m)}$ with $b^{(m)}\in V_m$.
Applying \eqref{weight-one-pairing} to each $b^{(m)}$ gives
\[
 \chi(u[1]v)\one
 =\sum_{m>0}(b^{(m)})_m u
 =\eta([b])=\lambda(v)\one,
\]
where the last equality is \eqref{projection} with $\Lz b=v$.
\end{proof}

We will also use the following consequence of simplicity. For a
simple vertex operator algebra of CFT type and $v\in V$,
\begin{equation}\label{vacuum-detection}
 \chi(Y(b,x)v)=0\ \text{for all } b\in V
 \quad\Longrightarrow\quad v=0 .
\end{equation}
Indeed, if $v\ne0$, the vectors $b_jv$ span $V$ by
\cite[Corollary~4.2]{DM}. Their vacuum coefficients cannot all
vanish, since $\chi(\one)=1$.

\begin{proof}[Proof of surjectivity in Theorem~\ref{main}]
Fix $D\in\Der(V)$, and let $\lambda$ be the functional given by
Proposition~\ref{functional}. Choose $u\in V_1$ as in
Lemma~\ref{realize}. We compare the two expressions for
$\lambda(a\circ b)$.

Since $u[0]=u_0$, the bracket commutator formula gives
\[
 [u[1],Y[a,z]]=Y[u[1]a,z]+zY[u_0a,z].
\]
Using $a\circ b=\Res_z\gamma(z)Y[a,z]b$, we obtain
\[
 u[1](a\circ b)
 =(u[1]a)\circ b+a\circ(u[1]b)
   +\Res_z z\gamma(z)Y[u_0a,z]b.
\]
The first two terms belong to $O(V)$ and are killed by $\chi$.
Thus Lemma~\ref{realize} and \eqref{lambda} give
\[
 \beta_{u_0}(a,b)
 =\chi(u[1](a\circ b))
 =\lambda(a\circ b)
 =-\beta_D(a,b).
\]
It follows that $F=D+u_0$ satisfies $\beta_F(a,b)=0$ for all
$a,b\in V$. It remains to show that $F=0$.

Let $a\in V_m$ and $b\in V_n$. Only $(Fa)_{m+n-1}b$ can have
weight zero. Taking the vacuum coefficient in the definition of
$\beta_F$ therefore gives
\[
 0=\left.\frac{d}{ds}\binom{s}{m+n+1}\right|_{s=m}
      \chi((Fa)_{m+n-1}b).
\]
The polynomial $\binom{s}{m+n+1}$ has a simple zero at $s=m$,
so $\chi((Fa)_{m+n-1}b)=0$. Hence $\chi(Y(Fa,x)b)=0$ for every
$b\in V$. Skew symmetry and $\chi T=0$ give
$\chi(Y(b,x)Fa)=\chi(Y(Fa,-x)b)=0$. By \eqref{vacuum-detection},
$Fa=0$. Thus $F=0$ and $D=(-u)_0$, proving surjectivity.
\end{proof}

\section{Faithfulness and finite automorphism groups}\label{group-sec}

The converse direction of the finiteness criterion requires that a nonzero vector in $V_1$ acts nontrivially on $V$. We give an argument using only rationality and the CFT grading.

\begin{lemma}\label{faithful}
Let $V$ be a rational vertex operator algebra of CFT type. If $h\in V_1$ and $h_0=0$ on $V$, then $h=0$.
\end{lemma}
\begin{proof}
Since $h[0]=h_0=0$, the bracket commutator formula gives
\[
 [h[1],Y[a,z]]=Y[h[1]a,z].
\]
Multiplying by $\phi(z)$ or $\gamma(z)$ and taking the residue
shows that $h[1]$ satisfies the derivation rule for both $*$ and
$\circ$. In particular, $h[1]$ preserves $O(V)$ and induces a
derivation $\delta$ of $A(V)$ given by
$\delta(\bar a)=\overline{h[1]a}$, where $\bar a=a+O(V)$.

Rationality implies that $A(V)$ is finite dimensional and semisimple
\cite[Theorem~2.2.3]{Zhu} (see also \cite[Theorem~4.10(a)]{DLM}
with $n=0$). The ideal $I=\ker\chi$ is therefore a
direct sum of matrix algebras and satisfies $I^2=I$. For $x,y\in I$,
both $\delta(x)y$ and $x\delta(y)$ lie in $I$, so
$\delta(I)=\delta(I^2)\subseteq I$. Since
$\bar a-\chi(\bar a)\bar\one\in I$ and $\delta(\bar\one)=0$, we get
\[
 \chi(h[1]a)=\chi(\delta(\bar a))=0\qquad(a\in V).
\]

For $a\in V_m$ with $m>0$, we have $a_jh\in V_{m-j}$. Hence
\eqref{weight-one-pairing} gives
\[
 \chi(Y(a,x)h)=x^{-m-1}\chi(a_mh)
              =m x^{-m-1}\chi(h[1]a)=0.
\]
For $a\in V_0=\C\one$, the same equality $\chi(Y(a,x)h)=0$
follows from $Y(\one,x)h=h\in V_1$. It therefore holds for every $a\in V$.

By rationality, the adjoint module $V$ is a direct sum of irreducible
submodules. Their weight zero spaces form a direct sum equal to
$\C\one$, so one submodule contains $\one$. It then contains every
$a=a_{-1}\one$ and hence equals $V$. Thus $V$ is simple, and
\eqref{vacuum-detection} gives $h=0$.
\end{proof}

\begin{proof}[Proof of Theorem~\ref{main}]
Surjectivity was proved in Section~\ref{realize-sec}. Injectivity is Lemma~\ref{faithful}. This proves part (1).

Part (2) was proved in Section~\ref{identity-sec}. Finally, let $G=\Aut(V)$, and let $\dim G$ denote the dimension of $G$
as a complex algebraic group. By Section~2 and part~(1), we have
\[
  \dim G=\dim\Lie(G)=\dim\Der(V)=\dim V_1.
\]
A complex algebraic group of finite type is finite if and only if it
has dimension zero. Hence $\Aut(V)$ is finite if and only if $V_1=0$. This proves part (3).
\end{proof}

\begin{example}
Let $L=\Z\alpha$ with $(\alpha,\alpha)=2$, and take $V=V_L$ with its usual conformal vector. This lattice vertex operator algebra is simple and of CFT type. It is regular \cite{DLMreg}, hence rational and $C_2$-cofinite \cite{ABD}. By \cite{DGR}, $\Aut(V_L)\cong\mathrm{PSL}_2(\C)$. This group is connected and infinite. Theorem~\ref{main}(2) therefore gives $\Inn(V_L)=\Aut(V_L)$, with trivial outer quotient.
\end{example}

\section*{Acknowledgments}
The author developed the idea of combining the two degree zero operators
$L(0)$ and $D$ to form the family $H_t=L(0)+tD$.
An AI system was used to search the literature, prepare the manuscript,
and help work out the details of the proofs.
The author is responsible for the mathematical content.

\end{document}